\documentstyle[12pt]{article}
\newcommand{\qed}{\hfill\rule{4pt}{8pt}\par\vspace{\baselineskip}}
\newtheorem{de}{Definition}[section]

\newtheorem{pr}[de]{Proposition}

\def\bea{\begin{eqnarray*}}\def\eea{\end{eqnarray*}}

\begin{document}
\title{Relating the Connes-Kreimer and Grossman-Larson Hopf algebras\\
built on rooted trees}

\author {Florin Panaite\\Institute of Mathematics of the Romanian Academy\\
P. O. Box 1-764, RO-70700 Bucharest, Romania\\
e-mail: fpanaite@stoilow.imar.ro}
\date{}
\maketitle
\section{Introduction}
${\;\;\;}$In \cite{k}, Dirk Kreimer discovered the striking fact that the 
process of renormalization in quantum field theory may be described, in a 
conceptual manner, by means of certain Hopf algebras (which depend  
on the chosen renormalization scheme). A toy model was studied in detail by 
Alain Connes and Dirk Kreimer in \cite{ck1}; the Hopf algebra which occurs, 
denoted by $H_R$, is the polynomial algebra in an infinity of 
indeterminates, one for each rooted tree, but with a non-cocommutative 
comultiplication. Some operators, denoted by $N$ and $L$, have been defined on 
$H_R$. The first one is the natural growth operator, which acts as a 
derivation; it defines some elements $\delta _k$, for $k\geq 1$, which 
provide the link between $H_R$ and another Hopf algebra, introduced in  
\cite{cm} by Alain Connes and Henri Moscovici in a completely different 
context, namely in noncommutative geometry. The operator $L$ is a solution 
of the ``Hochschild equation'', and the pair $(H_R, L)$ is characterized 
as the solution of a universal problem in Hochschild cohomology. It was 
proved also in \cite{ck1} that $H_R$ is in duality with the universal 
enveloping algebra  
of a certain Lie algebra $\cal L^1$, which has a linear basis indexed by all 
(non-empty)  
rooted trees. Let us note that this Lie algebra $\cal L^1$ appeared also, 
very recently, in \cite{cl}, in the context of pre-Lie algebras and the 
operad of rooted trees.\\
${\;\;\;}$In this note we would like to draw the attention to another Hopf 
algebra built on rooted trees, introduced ten years ago by Robert Grossman 
and Richard Larson in \cite{gl1} (see also their survey \cite{gl2}). 
This Hopf algebra (denoted by $A$ in what  
follows) has a $linear$ basis consisting of all (non-empty) 
rooted trees, a noncommutative 
product, and is a cocommutative graded connected Hopf algebra, hence, by 
the Milnor-Moore theorem, it is the universal enveloping algebra of 
the Lie algebra of its primitive elements, which may also be described 
explicitly: $P(A)$ has a linear  
basis consisting of all rooted trees whose root has exactly one child. Using  
these properties of $A$, Grossman and Larson gave a Hopf algebraic proof of 
the classical result of Cayley on the number of rooted trees. The 
construction of the Hopf algebra $A$ is motivated also by some ideas 
concerning differential operators and differential equations, in particular 
the Runge-Kutta method in numerical analysis (see \cite{gl2}). Let us note 
also that actually the construction of Grossman and Larson is slightly 
more general: they associate a Hopf algebra to any family of trees which 
satisfy a certain list of axioms. Among these families are: the family 
of all rooted trees (this is the one who gives the Hopf algebra $A$) and 
the families of all ordered, heap-ordered and respectively labelled  
rooted trees. The construction of the Hopf algebra associated to such a 
family is similar to that of $A$. \\  
${\;\;\;}$As noted in \cite{b}, \cite{ck2}, the Hopf algebra $H_R$ may also 
be related to the Runge-Kutta method and the Butcher group, so it is very 
likely that there is a relation between $H_R$ and $A$. As we shall see, this 
relation is best expressed by the fact that the Hopf algebras $A$ and 
$U(\cal L^1)$ are isomorphic, which in turn is a consequence of a Lie algebra 
isomorphism between $P(A)$ and $\cal L^1$. This isomorphism is given by 
sending a rooted tree $t\in P(A)$ to the rooted tree in $\cal L^1$ 
obtained by deleting the root of $t$.\\
${\;\;\;}$We believe that this relation between $A$ and $H_R$ may be 
useful for a better understading of both these Hopf algebras. On one hand, 
the advantage of the isomorphism between $A$ and $U(\cal L^1)$, which 
allows one to work with $A$ instead of $U(\cal L^1)$, is clear, since 
we know on $A$ a very explicit $linear$ basis, more manageable than the 
PBW basis of $U(\cal L^1)$. On the 
other hand, some known results for one of the Hopf algebras $A$ and $H_R$ 
may serve as a motivation and inspiration for obtaining similar results for 
the other. We shall make here a first step in this direction, by studying 
two natural operators on $A$. The first is the natural growth operator $N$ 
(defined exactly as the one introduced by Connes and Kreimer for $H_R$), 
which will turn out to be a $coderivation$ on $A$; the sequence 
$\{x_k\}_{k\geq 0}$    
defined by $N$ will turn out to generate a commutative cocommutative Hopf 
subalgebra of $A$, isomorphic to the polynomial algebra in one 
indeterminate, with its usual Hopf algebra structure. A nice feature of 
$N$ (considered on $A$) is that, for any rooted tree $t$, $N(t)$ may be 
described as the product (in $A$) between the rooted tree with two vertices 
and $t$.  
The second one, 
denoted by $M$, is in some sense dual to the operator $L$ on $H_R$: we 
shall prove that $M$ is a derivation (the right $A$-module structure on $A$ 
being the one induced by $\varepsilon $) and that the transpose of $M$ gives 
a solution of the Hochschild equation on the finite dual Hopf algebra 
$A^0$.      
\section{The relation between $H_R$ and $A$}
${\;\;\;}$Throughout, $k$ will be a fixed field of characteristic zero 
and all algebras, linear spaces etc. will be over $k$; unadorned 
$\otimes $ means $\otimes _k$.\\
${\;\;\;}$We start by recalling some facts from \cite{ck1}, \cite{gl1}, 
\cite{gl2}, to which we refer for the terminology and more details (the 
reader will find also in \cite{ck1} some nice pictures of rooted trees and 
the operations which may be performed with them).\\
${\;\;\;}$A $rooted$ $tree$ $t$ is a connected and simply-connected set of  
oriented edges and vertices such that there is exactly one distinguished 
vertex with no incoming edges, called the $root$ of $t$. Every edge connects 
two vertices. The $fertility$ $f(v)$ of a vertex $v$ is the number of edges 
outgoing from $v$, that is the number of $children$ of $v$. A $forest$ is a 
finite set of rooted trees.\\
${\;\;\;}$We denote by $B_-$ the operator which assigns to a rooted tree $t$ 
a forest, by removing the root of $t$, and by $B_+$ the operator which maps a 
forest consisting of $n$ rooted trees $t_1,...,t_n$ to a new rooted tree $t$  
which has a root $r$ with fertility $f(r)=n$ which connects to the $n$ 
roots of  
$t_1,...,t_n$. Obviously, $B_+(B_-(t))=B_-(B_+(t))=t$ for any rooted tree 
$t$. We also set $B_-(e)=\emptyset $, $B_+(\emptyset )=e$, where $e$ is the 
rooted tree with only one vertex and $\emptyset $ is the empty tree.\\
${\;\;\;}$If $t$ is a rooted tree, an $elementary$ $cut$ is a cut of $t$ at 
a single chosen edge, and an $admissible$ $cut$ is a set of elementary cuts 
such that any path from any vertex of $t$ to the root of $t$ contains at 
most one  
elementary cut. If $c$ is an admissible cut, we denote by $|c|$ the number of 
elementary cuts of $c$. If we perform an admissible cut $c$ in a rooted tree  
$t$, we obtain a forest, denoted by $P^c(t)$, consisting of the $cut$ 
$branches$ of $t$, and a $trunk$, denoted by $R^c(t)$, which is the branch 
which remains (it is the one which contains the root of $t$).\\
${\;\;\;}$We can now define the Connes-Kreimer Hopf algebra over $k$, denoted 
by $H_R$. As an algebra, it is the polynomial algebra in an infinity of 
indeterminates, one for each (non-empty) rooted tree (we denote also by 
$t$ the  
indeterminate corresponding to the rooted tree $t$). The unit is denoted 
by 1 (it corresponds to the empty tree). The comultiplication $\Delta $ is 
defined by:
$$\Delta (1)=1\otimes 1$$
$$\Delta (t)=1\otimes t+t\otimes 1+\sum _c P^c(t)\otimes R^c(t)$$
for any rooted tree $t$, where the sum is over all admissible cuts $c$ of $t$  
(with $|c|\geq 1$) and $P^c(t)$ is the monomial corresponding to the forest  
$P^c(t)$ (as a general rule, we identify any forest with its monomial). An  
alternative recursive description of $\Delta (t)$ is
$$\Delta (t)=t\otimes 1+(id\otimes B_+)(\Delta (B_-(t))$$
The counit is given by 
$$\varepsilon (1)=1$$
$$\varepsilon (t)=0$$
for any rooted tree $t$. The antipode is given iteratively by 
$$S(1)=1$$
$$S(t)=-t-\sum _cS(P^c(t))R^c(t)$$
for any rooted tree $t$, where the sum is over all admissible cuts $c$ of $t$  
(with $|c|\geq 1$).\\
${\;\;\;}$Define now the operator $N$ on $H_R$ (the natural growth operator) 
which maps a rooted tree $t$ with $n$ vertices to a sum $N(t)$ of $n$ rooted  
trees $t_i$, each having $n+1$ vertices, by attaching one more outgoing edge  
and vertex to each vertex of $t$ (the root remains the same under this 
operation). On products of rooted trees, $N$ acts, by definition, as a 
derivation. For any $k\geq 1$, define elements $\delta _k\in H_R$ by 
$\delta _1=e$, $\delta _{k+1}=N(\delta _k)$, that is $\delta _{k+1}=
N^k(e)$ for all $k\geq 1$. We recall that we have denoted by $e$ the tree 
with one vertex (in order to be consistent with \cite{gl1}; in \cite{ck1} by 
$e$ is denoted the unit of $H_R$). These elements are very important in 
\cite{ck1}, because they provide the link between $H_R$ and the 
Connes-Moscovici Hopf algebra introduced in \cite{cm}. They generate a 
(non-cocommutative) Hopf subalgebra of $H_R$.\\
${\;\;\;}$Define now the operator $L:H_R\rightarrow H_R$ as the unique  
linear map satisfying the condition $L(t_1...t_m)=t$ for any rooted trees 
$t_1,...,t_m$, where $t$ is the rooted tree obtained by connecting a new root 
to the roots of $t_1,...,t_m$. Obviously it agrees with the map $B_+$ 
introduced above. It was shown in \cite{ck1} that the operator $L$ satisfies 
the so-called ``Hochschild equation''
$$\Delta \circ L=L\otimes 1+(id\otimes L)\circ \Delta $$
(this is a 1-cocycle condition) and the pair $(H_R, L)$ has the 
following universal property: if $(H_1, L_1)$ is a pair with $H_1$ a 
commutative Hopf algebra and $L_1:H_1\rightarrow H_1$ a linear map satisfying  
the Hochschild equation on $H_1$, then there exists a unique Hopf algebra map 
$\rho :H_R\rightarrow H_1$ such that $L_1\circ \rho =\rho \circ L$.\\ 
${\;\;\;}$Let $\cal L^1$ be the linear span of the elements $Z_t$, indexed by 
all (non-empty) rooted trees. Define an operation on $\cal L^1$ by  
$$Z_{t_1}*Z_{t_2}=\sum _t n(t_1,t_2;t)Z_t$$
where the integer $n(t_1,t_2;t)$ is the number of admissible cuts $c$ of $t$ 
with $|c|=1$ such that the cut branch is $t_1$ and the trunk is $t_2$ (note 
that this operation is $not$ associative).  
Define then a bracket on $\cal L^1$ by 
$$[Z_{t_1}, Z_{t_2}]=Z_{t_1}*Z_{t_2}-Z_{t_2}*Z_{t_1}$$
Then it was proved in \cite{ck1} that $(\cal L^1, [,])$ is a Lie algebra and 
moreover there is a Hopf duality between $H_R$ and the universal enveloping  
algebra of $\cal L^1$.\\[2mm]
${\;\;\;}$We recall now from \cite{gl1}, \cite{gl2} the structure of the 
Grossman-Larson Hopf algebra, which will be denoted in what follows by $A$. 
It has a linear basis consisting of all (non-empty) rooted trees. The unit  
is the tree $e$ with only one vertex. The multiplication on the basis is 
given as follows: let $t_1$ and $t_2$ be two rooted trees, let $s_1,...,s_r$ 
be the children of the root of $t_1$, let $n$ be the number of vertices of  
$t_2$; then there are $n^r$ ways to attach the $r$ subtrees of $t_1$ which  
have $s_1,...,s_r$ as roots to the tree $t_2$ by making each $s_i$ the child 
of some vertex of $t_2$. The product $t_1t_2$ is defined to be the sum 
of these $n^r$ rooted trees (note that this product is $not$ commutative).\\
${\;\;\;}$The coalgebra structure of $A$ is given as follows. If $t$ is a 
rooted tree whose root has children $s_1,...,s_r$, the coproduct 
$\Delta (t)$ is the sum of the $2^r$ terms $t_1\otimes t_2$, where the 
children of the root of $t_1$ and the children of the root of $t_2$ range 
over all $2^r$ possible partitions of the children of the root of $t$ into 
two subsets. If $t=e$, then $\Delta (e)=e\otimes e$. The counit 
$\varepsilon $ is given by $\varepsilon (e)=1$, $\varepsilon (t)=0$ if 
$t\neq e$. Obviously $\Delta $ is cocommutative. Moreover, $A$ is a graded 
connected bialgebra, the component $A_n$ of degree $n$ having as basis all  
trees with $n+1$ vertices. Being a graded connected cocommutative bialgebra,  
A is a Hopf algebra and by the Milnor-Moore theorem $A$ is the universal 
enveloping algebra of its primitives, $A\simeq U(P(A))$, where  
$P(A)=\{a\in A/\Delta (a)=e\otimes a+a\otimes e\}$. There is an explicit  
description of $P(A)$: it has a basis consisting of all rooted trees whose 
root has exactly one child.\\[3mm]
${\;\;\;}$We can state now the result which expresses the relation between  
the Hopf algebras $H_R$ and $A$. 
\begin{pr} The Lie algebras $\cal L^1$ and $P(A)$ are isomorphic, hence $A$ 
is isomorphic to $U(\cal L^1)$ as Hopf algebras.
\end{pr}
{\bf Proof:} Define $\varphi :P(A)\rightarrow \cal L^1$ as the unique linear 
map which on the basis of $P(A)$ acts as follows: if $t\in P(A)$ is a rooted 
tree, then $\varphi (t)=Z_{B_-(t)}$. Recall that $B_-(t)$ is the rooted tree  
obtained by deleting the root of $t$ (here it is a tree,  
since the root of $t$ has exactly one child). Obviously, $\varphi $ is a 
linear isomorphism, its inverse being the map $\psi :\cal L^1\rightarrow 
P(A)$, $\psi (Z_t)=B_+(t)$ for any rooted tree $t$. \\
It remains to prove that $\varphi $ is a Lie algebra map. Let 
$t_1,t_2\in P(A)$ be two rooted trees. In $P(A)$, we have $[t_1,t_2]=t_1t_2-
t_2t_1$. By the definition of the multiplication of $A$, we obtain that 
$$t_1t_2=B(t_1t_2) +\sum _it^i$$
$$t_2t_1=B(t_2t_1)+\sum _jT^j$$
where $B(t_1t_2)$ is the rooted tree obtained by identifying the roots of 
$t_1$ and $t_2$, and each $t^i$ is obtained by identifying the root of $t_1$ 
with a vertex of $t_2$, except for the root of $t_2$ (and similarly for 
$t_2t_1$), so all the rooted trees $t^i$ and $T^j$ are in $P(A)$. 
Obviously $B(t_1t_2)=B(t_2t_1)$, hence 
$$[t_1,t_2]=\sum _it^i-\sum _jT^j$$
We obtain that $\varphi ([t_1,t_2])=\sum _i\varphi (t^i)-
\sum _j\varphi (T^j)$. From the definition of the operation $*$ on 
$\cal L^1$, it is easy to see that
$$\sum _i\varphi (t^i)=\varphi (t_1)*\varphi (t_2)$$
$$\sum _j\varphi (T^j)=\varphi (t_2)*\varphi (t_1)$$
hence we obtain
$$\varphi ([t_1,t_2])=[\varphi (t_1), \varphi (t_2)]$$
that is $\varphi $ is a Lie algebra map. \qed

${\;\;\;}$Define now the natural growth operator $N$ on $A$, by the same 
formula as the one defined by Connes and Kreimer on $H_R$, that is $N$ is the 
linear map $N:A\rightarrow A$ such that, for any rooted tree $t$, $N(t)$ is 
the sum of the rooted trees obtained from $t$ by attaching one more 
outgoing edge and vertex to each vertex of $t$. The properties of $N$ are 
collected in the following
\begin{pr} $(1)\; N(e)b=N(b)$ for all $b\in A$.\\[2mm]
$(2)\; N(ab)=N(a)b$ for all $a,b\in A$. \\[2mm]
$(3)\; N^k(e)b=N^k(b)$ for all $k\geq 1, b\in A$.\\[2mm]
$(4)\; N^k(e)N(e)=N^{k+1}(e)=N(e)N^k(e)$ for all $k\geq 1$.\\[2mm]
$(5)\; N$ is a coderivation, that is, for all $b\in A$, we have 
$$\Delta (N(b))=(id\otimes N+N\otimes id)(\Delta (b))$$
\end{pr}
{\bf Proof:} $(1)$ obviously, for any rooted tree $t$, we have $N(e)t=N(t)$, 
from which we obtain $N(e)b=N(b)$ for all $b\in A$. \\
$(2)$ using $(1)$, we have $N(ab)=N(e)ab=N(a)b$.\\
$(3)$ follows easily by induction, and $(4)$ follows from $(3)$.\\
$(5)$ by $(1)$, we obtain that $\Delta (N(b))=\Delta (N(e)b)=  
\Delta (N(e))\Delta (b)$. We shall use the $\Sigma $-notation, so we write 
$\Delta (b)=\sum b_{(1)}\otimes b_{(2)}$; since $N(e)$ is a primitive 
element, we have then:
$$\Delta (N(b))=(e\otimes N(e)+N(e)\otimes e)(\sum b_{(1)}\otimes b_{(2)})=$$
$$\sum b_{(1)}\otimes N(e)b_{(2)}+\sum N(e)b_{(1)}\otimes b_{(2)}=$$
$$\sum b_{(1)}\otimes N(b_{(2)})+\sum N(b_{(1)})\otimes b_{(2)}=$$
$$(id\otimes N+N\otimes id)(\Delta (b))$$\qed

${\;\;\;}$Now, for any $k\geq 0$, define the element $x_k=N^k(e)\in A$. These 
elements are analogous to the elements $\delta _k$ of \cite{ck1}.
\begin{pr} The elements $x_k$ have the following properties: \\[2mm]
$(1)\; x_mx_n=x_nx_m=x_{m+n}$\\[2mm]
$(2)\; x_0=e$\\[2mm]
$(3)\; \varepsilon (x_m)=\delta _{0,m}$\\[2mm]
$(4)\; \Delta (x_m)=\sum _{i=0}^m \left (\begin{array}{c} 
m\\i \end{array}\right) x_i\otimes x_{m-i}$\\[2mm]
$(5)\; S(x_m)=(-1)^mx_m$ \\[2mm]
for all $m,n\geq 0$, where $S$ is the antipode of $A$. Moreover, the 
elements $\{x_k\}_{k\geq 0}$ are linearly independent.    
Hence, the subspace of $A$ generated by the elements $x_k$ with $k\geq 0$ is 
a commutative cocommutative Hopf subalgebra of $A$, isomorphic (as a Hopf 
algebra) to the polynomial algebra $k[X]$ with its usual Hopf algebra 
structure.  
\end{pr}
{\bf Proof:} $(1)$ follows easily from the previous proposition; $(2)$ and 
$(3)$ are obvious; $(4)$ and $(5)$ follow by induction, using the facts 
that $x_{m+1}=N(e)x_m$ and $N(e)$ is primitive. Since $x_k\in A_k$ for all 
$k\geq 0$ (hence $deg (x_m)\neq deg (x_n)$ if $m\neq n$) it follows that the  
elements $\{x_k\}_{k\geq 0}$ are linearly independent.  
The isomorphism between 
the Hopf subalgebra given by the elements $x_k$ and $k[X]$ is determined by 
$x_k\mapsto X^k$ for all $k\geq 0$.\qed

${\;\;\;}$Define now the $k$-linear map $M:A\rightarrow A$, by $M(e)=0$ and 
$M(t)=tN(e)$ for all rooted trees $t\neq e$. As we shall see, this operator  
$M$ is in some sense dual to the operator $L$ on $H_R$. 
\begin{pr} $(1)\; M(b)=(b-\varepsilon (b))N(e)$ for all $b\in A$.\\[2mm]
$(2)\; M(ab)=aM(b)+M(a)\varepsilon (b)$ for all $a,b\in A$\\[2mm]
(that is, $M$ is a derivation from $A$ into the bimodule $A$, where the left 
$A$-module structure of $A$ is the one given by multiplication and the right 
$A$-module structure is given via $\varepsilon $).\\[2mm]
$(3)\; \Delta (M(b))=(id \otimes M+M\otimes id)(\Delta (b))+
(b-\varepsilon (b))\otimes N(e)+N(e)\otimes (b-\varepsilon (b))$\\[2mm]
for all $b\in A$.
\end{pr}
{\bf Proof:} $(1)$ Let $b\in A$, written as $b=b_0e+\sum _{t_i\neq e}b_it_i$, 
with $b_0, b_i\in k$. Then we have: 
$$M(b)=b_0M(e)+\sum b_iM(t_i)=\sum b_it_iN(e)=$$
$$\sum b_it_iN(e)+b_0N(e)-b_0N(e)=bN(e)-b_0N(e)=(b-\varepsilon (b))N(e)$$
$(2)$ using $(1)$, we have $M(ab)=(ab-\varepsilon (a)\varepsilon (b))N(e)$, 
and 
$$aM(b)+M(a)\varepsilon (b)=a(b-\varepsilon (b))N(e)+(a-\varepsilon (a))
\varepsilon (b)N(e)=$$
$$(ab-\varepsilon (b)a+\varepsilon (b)a-\varepsilon (a)\varepsilon (b))N(e)=
(ab-\varepsilon (a)\varepsilon (b))N(e)$$
$(3)$ write $\Delta (b)=\sum b_{(1)}\otimes b_{(2)}$; then we have:
$$\Delta (M(b))=\Delta ((b-\varepsilon (b))N(e))=\Delta (b)\Delta (N(e))-
\varepsilon (b)\Delta (N(e))=$$
$$\sum b_{(1)}N(e)\otimes b_{(2)}+\sum b_{(1)}\otimes b_{(2)}N(e)-
\varepsilon (b)\otimes N(e)-N(e)\otimes \varepsilon (b)=$$
$$\sum M(b_{(1)})\otimes b_{(2)}+\sum \varepsilon (b_{(1)})N(e)
\otimes b_{(2)}+\sum b_{(1)}\otimes M(b_{(2)})+$$
$$+\sum b_{(1)}\otimes 
\varepsilon (b_{(2)})N(e)-\varepsilon (b)\otimes N(e)-N(e)\otimes 
\varepsilon (b)=$$
$$(id\otimes M+M\otimes id)(\Delta (b))+(b-\varepsilon (b))\otimes N(e)+
N(e)\otimes (b-\varepsilon (b))$$\qed   
${\;\;\;}$Recall from \cite{sw}, \cite{m} that we can associate to 
$A$ the so-called  
$finite\;dual$, which is also a Hopf algebra, and which consists of the 
elements $f\in A^*$ such that $Ker (f)$ contains a cofinite ideal of $A$. 
If we denote by $m$ the multiplication of $A$, then an element $f\in A^*$ 
belongs to $A^0$ if and only if $m^*(f)\in A^*\otimes A^*$, which in turn 
is equivalent to the fact that there exist some elements $f_i, f'_i\in A^*$,  
with $i$ in some finite set, such that $f(ab)=\sum f_i(a)f'_i(b)$ for all 
$a,b\in A$. Moreover, we have that $m^*(A^0)\subseteq A^0\otimes A^0$ and 
the comultiplication of $A^0$ is $\Delta =m^*|A^0$.\\
${\;\;\;}$Now, let $f\in A^0$, $m^* (f)=\sum f_i\otimes f'_i$, with 
$f_i,f'_i\in A^0$. By  
using the condition $M(ab)=aM(b)+M(a)\varepsilon (b)$ satisfied by $M$, 
we can compute:
$$f(M(ab))=f(M(a))\varepsilon (b)+f(aM(b))$$
for all $a,b\in A$, which may be rewritten as 
$$f(M(ab))=f(M(a))\varepsilon (b)+\sum f_i(a)f'_i(M(b))$$
that is $$M^*(f)(ab)=M^*(f)(a)\varepsilon (b)+\sum f_i(a)M^*(f'_i)(b)$$
hence $M^*(f)\in A^0$. So, we have $M^*(A^0)\subseteq A^0$, and we denote 
by $M^0$ the restriction of $M^*$ to $A^0$. Also, for $f\in A^0$, since 
$M^*(f), \varepsilon, f_i, M^*(f'_i)\in A^0$ for all $i$, we obtain 
finally that in $A^0$ the following relation holds:
$$\Delta (M^0(f))=M^0(f)\otimes \varepsilon +(id \otimes M^0)(\Delta (f))$$
which expresses the fact that $M^0$ is a solution for the Hochschild  
equation on $A^0$. Hence, since $A^0$ is commutative (because $A$ is 
cocommutative),  
by the universal property of the pair $(H_R, L)$ 
we obtain the following 
\begin{pr} There exists a unique Hopf algebra map $\rho :H_R\rightarrow A^0$ 
such that $M^0\circ \rho=\rho \circ L$.
\end{pr}                        

\end{document}